\documentclass[12pt]{amsart}
\usepackage{a4wide}
\usepackage{amssymb}
\usepackage{amsthm}
\usepackage{amsmath}
\usepackage{amscd}
\usepackage{verbatim}
\usepackage[all]{xy}
\usepackage{longtable, lscape}

\theoremstyle{plain}
\newtheorem{theorem}{Theorem}[section]
\newtheorem{corollary}[theorem]{Corollary}
\newtheorem{lemma}[theorem]{Lemma}
\newtheorem{proposition}[theorem]{Proposition}

\newtheorem{conjecture}[theorem]{Conjecture}
\theoremstyle{definition}
\newtheorem{definition}[theorem]{Definition}

\theoremstyle{remark}
\newtheorem{remark}{Remark}

\newcommand{\R}{\mathbb{R}}
\newcommand{\Q}{\mathbb{Q}}
\newcommand{\Z}{\mathbb{Z}}

\newcommand{\C}{\mathbb{C}}

\renewcommand{\H}{\mathbb{H}}
\renewcommand{\P}{\mathbb{P}}

\newcommand{\T}{\mathbb{T}}

\newcommand{\leg}[2]{\left( \frac{#1}{#2} \right)}
\newcommand{\kzxz}[4]{\left(\begin{smallmatrix} #1 & #2 \\ #3 & #4\end{smallmatrix}\right) }
\newcommand{\kabcd}{\kzxz{a}{b}{c}{d}}

\newcommand{\calD}{\mathcal{D}}

\newcommand{\calO}{\mathcal{O}}
\newcommand{\calP}{\mathcal{P}}
\newcommand{\calQ}{\mathcal{Q}}

\newcommand{\calW}{\mathcal{W}}

\newcommand{\frake}{\mathfrak e}

\newcommand{\eps}{\varepsilon}

\newcommand{\sgn}{\operatorname{sgn}}

\newcommand{\Sl}{\operatorname{SL}}

\newcommand{\Mp}{\operatorname{Mp}}
\newcommand{\Orth}{\operatorname{O}}

\newcommand{\Aut}{\operatorname{Aut}}

\newcommand{\dv}{\operatorname{div}}

\newcommand{\res}{\operatorname{res}}

\newcommand{\ccl}{\mathcal{CL}}

\newcommand{\Div}{\operatorname{Div}}

\begin{document}

\title{Harmonic Maass forms and periods}

\date{\today}
\thanks{The author is partially supported by DFG grant BR-2163/2-1.}

\author[Jan H. Bruinier]{Jan Hendrik Bruinier}


\address{Fachbereich Mathematik,
Technische Universit\"at Darmstadt, Schlossgartenstrasse 7, D--64289
Darmstadt, Germany} \email{bruinier@mathematik.tu-darmstadt.de}

\begin{abstract}
  According to Waldspurger's theorem, the coefficients of
  half-integral weight eigenforms are given by central critical values
  of twisted Hecke $L$-functions, and therefore by periods.  Here we
  prove that the coefficients of weight $1/2$ harmonic Maass forms are
  determined by periods of algebraic differentials of the third kind on
  modular and elliptic curves.
\end{abstract}

\maketitle

\section{Introduction}
\label{sect:intro}


The Fourier expansions of half integral weight modular forms serve as
generating series of important number theoretic functions, such
as representation numbers of quadratic forms and
class numbers of imaginary quadratic fields.
%
The Shimura correspondence \cite{Sh} provides a map from holomorphic
modular forms of half-integral weight $k+1/2$ to forms of weight $2k$,
which is compatible with the action of the Hecke algebra.  A
celebrated result of Waldspurger \cite{Wa} and Kohnen--Zagier
\cite{KZ} says that the coefficients of square-free index of a Hecke
eigenform $g$ of weight $k+1/2$ are given by the central critical
values of the quadratic twists of the Hecke $L$-function of the
Shimura lift of $g$.  There are many applications, in particular in
connection with the Birch and Swinnerton-Dyer conjecture.

Variants of Waldspurger's theorem also hold for non-holomorphic
modular forms. Katok and Sarnak considered the Shimura correspondence
between Maass cusp forms of weights $1/2$ and $0$. They showed that the
coefficients of a Maass eigenform $\varphi$ of weight $1/2$ are determined
by cycle integrals and CM values of the Maass eigenform of weight $0$
corresponding to $\varphi$ under the Shimura lift. Such results are crucial
in Duke's work on the equidistribution of CM
points on modular curves \cite{Du}.

In the present paper we study the coefficients of harmonic weak Maass forms
of weight $1/2$.  We prove that they are given by period
integrals of algebraic differential forms on modular and elliptic
curves. This leads to a refinement and a strengthening of recent results
of Ono and the author (see \cite{BO}, and Theorem \ref{Lval2} here) 
relating the coefficients to central values and
derivatives of twisted Hecke $L$-functions.

We now describe the content of the present paper in more detail.
Throughout, for $\tau$ in the complex upper half plane $\H$, we let $\tau=u+iv$, where $u, v\in \R$,
and we let $q:=e^{2 \pi i \tau}$. A
harmonic weak\footnote{For brevity we will often drop the attribute
  ``weak'' in this paper.} Maass form of weight $k\in \frac{1}{2}\Z$
on $\Gamma_0(N)$ (with $4\mid N$ if $k\in \frac{1}{2}\Z\setminus \Z$)
is a smooth function on $\H$, the upper half of the complex plane,
which satisfies:
\begin{enumerate}
\item[(i)]
 $f\mid_k\gamma = f$ for all $\gamma\in \Gamma_0(N)$;
\item[(ii)] $\Delta_k f =0 $, where $\Delta_k$ is the weight $k$
hyperbolic Laplacian on $\H$;
\item[(iii)]
there is a polynomial $P_f(q)=\sum_{n\leq 0} c^+(n)q^n \in \C[q^{-1}]$
such that $f(\tau)-P_f(q) \to 0$ as $v\to\infty$.
Analogous conditions are required at all
cusps (see \cite{BF}).
\end{enumerate}
The polynomial $P_f$ is called the
principal part of $f$ (at the cusp $\infty$).

Such a harmonic Maass form $f$ has a Fourier expansion of the form
\begin{align}
\label{intro:ffourier}
f(\tau)=\sum_{n\geq n_0} c^+(n)q^n + \sum_{n<0} c^-(n)\Gamma(1-k,4\pi |n| v)q^n,
\end{align}
where $\Gamma(a,x)$ denotes the incomplete gamma function.  The
differential operator $\xi(f)=2iv^k\overline{\frac{\partial
    f}{\partial \bar \tau}}$ takes harmonic Maass forms of weight $k$
to cusp forms of dual weight $2-k$. 
The coefficients $c^-(n)$ of the non-holomorphic part of $f$ are
up to an elementary factor equal to the coefficients of the cusp form
$\xi(f)$, while
the coefficients $c^+(n)$ of the holomorphic part are rather
mysterious.

Harmonic Maass forms have been a source of recent interest due to
their connection to Ramanujan's mock theta functions, see
e.g.~\cite{BriO}, \cite{BOR}, \cite{On}, \cite{Za}, \cite{Z1},
\cite{Z2}.  The mock theta functions correspond to harmonic Maass
forms of weight $1/2$ whose image under the differential operator
$\xi$ is a linear combination of unary weight $3/2$ theta
series. Their coefficients often have combinatorial interpretations.

In the present paper we consider (similarly as in \cite{BO}) those harmonic
Maass forms $f$ of weight $1/2$ for which $\xi(f)$ is orthogonal to
all unary theta series of weight $3/2$.  From an arithmetic
perspective, this case is of particular interest, since the Shimura
lift of $\xi(f)$ is a weight $2$ {\em cusp form}, which leads to a
connection to elliptic curves.



Our main
result is an explicit formula for the coefficients $c^+(n)$ in terms
of periods of algebraic differentials.
To illustrate it, we consider as an example the unique harmonic Maass
form $f_3$ of weight $1/2$ in the Kohnen plus space for
$\Gamma_0(4\cdot 37)$ whose principal part at the cusp $\infty$ is
given by $q^{-3}$. The coefficients $c^\pm(n)$ can be numerically
computed (see \cite{BS}), the first few are listed in
Table~\ref{table:1}.

\begin{table}[h]\begin{center}
\caption{\label{table:1} The first few coefficients of $f_3$}
\begin{tabular}{rcc}
\rule[-3mm]{0mm}{6mm}
$\Delta$ && $c^+(\Delta)$\\
\hline
\rule[-1mm]{0mm}{6mm}
$1$   &&  $-0.281761784989599568797560755375154934\dots$
      \\
\rule[-1mm]{0mm}{6mm}
$12$  &&  $-0.488527238262012252282270296073370716\dots$
      \\
\rule[-1mm]{0mm}{6mm}
$21$  &&  $-0.172739257232652756520820073970689922\dots$
      \\
\rule[-1mm]{0mm}{6mm}
$28$  &&  $\phantom{-} 0.678193995303947798284505784006694209\dots$
      \\
\rule[-1mm]{0mm}{6mm}
$33$  &&  $\phantom{-} 0.566302320159069981682205456692456226\dots$
\end{tabular}
\end{center}\end{table}

The image of $f_3$ under the differential operator $\xi$ is the cusp
form of weight $3/2$ whose Shimura lift is the eigenform $G$ of weight
$2$ and level 37 corresponding to the elliptic curve
\begin{align*}
E: y^2= 4x^3-4x+1
\end{align*}
of conductor $37$.
The Mordell-Weil group of $E$ is infinite cyclic
with generator $(0,-1)$.   The meromorphic differential
$\frac{-1}{x}\frac{dx}{y}$ on $E$ is of the third kind. It is regular
up to first order poles at the points $(0,-1)$ and $(0,1)$ with
residues $1$ and $-1$, respectively. It turns out that with the two periods
\begin{align*} 
\Theta&=\Re\left(\int_{E(\R)} \frac{-1}{2x}\, \frac{dx}{y}\right)=-1.68688450290973441728\dots,\\
\Omega&=\int_{E(\R)} \frac{dx}{y}=5.9869172924639192\dots,
\end{align*}
we have $c^+(1)=\Theta/\Omega$!

This identity is a consequence of a more general result on
differentials of the third kind on modular curves.  Let $G$
be a newform of weight $2$ for $\Gamma_0(N)$. For any Hecke operator
$T$, let $\lambda_G(T)$ be the corresponding eigenvalue. Throughout
this introduction we assume for simplicity that the level $N$ is a
prime and that $G$ is defined over $\Q$ and invariant under the Fricke
involution. Then the $L$-function of $G$ has an odd functional
equation. The case of general level, Fricke eigenvalue and field of
definition is treated in the body of the paper by working with vector
valued modular forms.

According to \cite[Lemma 7.3]{BO}, there exists a harmonic Maass form
$f$ of weight $1/2$ in the Kohnen plus space for $\Gamma_0(4N)$ such
that
\begin{itemize}
\item[(i)] $P_f\in \Z[q^{-1}]$ and and the constant term vanishes,
\item[(ii)] the Shimura lift of $\xi(f)$ is equal to $G$.
\end{itemize}
As before we denote the Fourier coefficients of $f$ by $c^\pm(n)$.
Waldspurger's theorem and the action of $\xi$ on the Fourier expansion imply that
for negative fundamental discriminants $\Delta$ with
$\leg{\Delta}{N}=1$ we have
$$
L(G,\chi_{\Delta},1)=C \sqrt{|\Delta|}\cdot
c^{-}(\Delta)^2,
$$
where $C>0$ is a constant which does not depend on $\Delta$.

To describe the coefficients $c^+(\Delta)$ of the holomorphic part we
consider differentials of the third kind associated with Heegner
divisors.  Let $d<0$ and $\Delta>0$ be fundamental discriminants which
are both squares modulo $N$. Let $\calQ_{d,N}$ be the set of
discriminant $d=b^2-4ac$ integral binary quadratic forms
$aX^2+bXY+cY^2$ with the property that $N\mid a$. Recall that there is
a twisted Heegner divisor $Z_{\Delta}(d)$ on the modular curve $X=X_0(N)$ defined by
\begin{align}
Z_\Delta(d) = \sum_{Q\in \calQ_{\Delta d,N}/\Gamma_0(N)}
\chi_\Delta (Q)\cdot  \frac{\alpha_Q}{w_Q}.
\end{align}
Here $\chi_\Delta$ denotes the generalized genus character
corresponding to the decomposition $\Delta\cdot d$ as in \cite{GKZ},
 $\alpha_Q$ is the unique root of $Q(x,1)$ in $\H$, and $w_Q$ denotes
 the order of the stabilizer  of $Q$ in $\Gamma_0(N)$.
The field of definition of $Z_\Delta(d)$ is 
$\Q(\sqrt{\Delta})$.

We associate a divisor $Z_\Delta(f)$ with the harmonic Maass form $f$ by putting
\begin{align}
Z_{\Delta}(f)=\sum_{n<0} c^+(n) Z_{\Delta}(n)\in
\Div^0(X).
\end{align}
According to \cite[Theorem 7.5]{BO}, the class of $Z_{\Delta}(f)$
lies in the $G$-isotypical component of the the Jacobian of $X$.

Recall that a differential of the third kind on $X$ is a
meromorphic $1$-form $\psi$ that has at most simple poles with
integral residues.  If $\psi$ has its poles at the points $P_j\in X$
with residues $c_j$, then the divisor
\[
\res(\psi)=\sum_{j} c_j\cdot P_j
\]
is called the residue divisor of $\psi$.  The Hecke algebra acts on
differentials of the third kind.  In Proposition \ref{prop:3} we will
show that there exists a unique differential of the third kind
$\zeta_\Delta(f)$ on $X$ with the following properties:
\begin{itemize}
\item[(i)]
The residue divisor is given by
$\res(\zeta_\Delta(f))=Z_\Delta(f)$.
\item[(ii)] The differential $\zeta_\Delta(f)$ is $G$-isotypical, that
  is, for all Hecke operators $T$, the differential
  $T\psi_\Delta(f)-\lambda_G(T)\psi_\Delta(f)$ is equal to
  $\frac{dF}{F}$ for a rational function
  $F\in \C(X)^\times$.
\item[(iii)] The first Fourier coefficient of $\zeta_\Delta(f)$ vanishes.
\end{itemize}
The differential $\zeta_\Delta(f)$  is called the normalized differential of the third kind for $Z_\Delta(f)$.
It is defined over $\Q(\sqrt{\Delta})$.

Let $H_1^+(X,\R)$ be the subspace of the first homology of $X$ which
is invariant under the involution induced by complex conjugation. The
Hecke algebra also acts on this space and decomposes it into
isotypical components.  Let $C_G$ be a generator of the $G$-isotypical
component. Our first main result (see Theorem \ref{thm:comp}) states:

\begin{theorem}
Let $\omega_G=2\pi i \cdot G(z)\,dz$, and let $\zeta_\Delta(f)$ be as before.
For any fundamental discriminant $\Delta>0$ which is a square modulo $N$ we have
\[
c^+(\Delta)= \frac{\Re\int_{C_G} \zeta_\Delta(f)}{\sqrt{\Delta}\int_{C_G} \omega_G}.
\]
\end{theorem}

\begin{remark}
Note that a result of Waldschmidt on the transcendence of periods of
differentials of the third kind (see Section 5.2 of \cite{W}, and
Theorem 2 of \cite{Sch}) implies that the right hand side is algebraic
if and only if $Z_\Delta(f)$ defines a  torsion point in the Jacobian.  Combining
this with the Gross-Zagier formula \cite{GZ}, it can be deduced  that
\[
c^+(\Delta)\in \Q\;\Leftrightarrow \; c^+(\Delta)\in \bar \Q \;\Leftrightarrow \; L'(G,\chi_\Delta,1)=0,
\]
which can be viewed as a weak version of Waldspurger's theorem for
central derivatives of quadratic twists of Hecke $L$-functions (see \cite{BO}).
\end{remark}

Let $E$ be an elliptic curve over $\Q$ corresponding to the newform
$G$.  To relate the coefficients of $f$ to periods of differentials on $E$ we
use a modular parameterization $\phi:X\to E$.  The divisor
$Z_\Delta(f)$ gives rise to a point $P_{\Delta}(f)=(x_\Delta,y_\Delta)$ in
$E(\Q(\sqrt{\Delta}))$.  It is mapped to its negative under the
non-trivial automorphism $\iota$ of $\Q(\sqrt{\Delta})$ and therefore defines
a rational point on the $\Delta$-th quadratic twist of $E$.
The differential
\begin{align}
\label{intro:eq:beta}
\beta_\Delta(f)=\frac{y_\Delta}{2(x-x_\Delta)} \cdot \frac{dx}{y}
\end{align}
on $E$ is of the third kind and has the residue divisor $\frac{1}{2}
(P_{\Delta}(f))-\frac{1}{2}(-P_{\Delta}(f))$. It is defined over
$\Q(\sqrt{\Delta})$ and $\iota$ takes it to its negative. We show
that its period over $E( \R)$ is related to the coefficient
$c^+(\Delta)$ (see Theorem \ref{thm:ec}):

\begin{theorem}
\label{intr:thm:ec}
The quantity
\begin{align}
\label{intro:diff}
 c^{+}(\Delta) - \frac{c_E }{ \sqrt{\Delta}\Omega(E) }\cdot \Re\int_{E(\R)}\beta_{\Delta}(f).
\end{align}
is rational. Here $\Omega(E)$ denotes the real period, and $c_E$ denotes the Manin constant of $E$.
\end{theorem}

Note that $\beta_\Delta(f)$ is determined by its residue divisor and the action of $\Aut(\C)$  only up to addition of a rational multiple of $\sqrt{\Delta}\frac{dx}{y}$.

It might be possible to refine Theorem~\ref{intr:thm:ec} by
working with a minimal Weierstrass model $\calW$ of $E$ over $\Z$. We
show that there exists a rational section
of the
sheaf $\Omega_{\calW/\Z}^1$ of relative differentials whose residue divisor is the Zariski
closure of $\frac{1}{2} (P_{\Delta}(f))-\frac{1}{2}(-P_{\Delta}(f))$
up to possible vertical components with bounded multiplicities above
$2$. We conjecture that if one replaces $\beta_\Delta(f)$ by such a
differential, then the analogue of \eqref{intro:diff} will be
contained in $\frac{1}{4}\Z$, see Conjecture \ref{con:ec}. We present
some numerical data supporting this in Section \ref{subsect:ex}.

The present paper is organized as follows: In Section
\ref{sect:difftk} we study differentials of the third kind on modular
curves.
In particular we consider their pairing
with the first homology and the action of the Hecke
algebra. In Section 3 we collect some facts on harmonic Maass forms and Heegner divisors.
We prove our main result by comparing the canonical differentials for certain Heegner divisors constructed in \cite{BO} with the differentials $\zeta_\Delta(f)$.
In Section 4 we consider the applications to rational elliptic curves.

The idea for the comparison of differentials in Section~3
is motivated by \cite[Section~3.5]{KonZ}.
I would like to thank Don Zagier for drawing my attention to this work.
I also thank him and Ken Ono for inspiring and helpful conversations.

\section{Differentials of the third kind}

\label{sect:difftk}

We
begin by recalling some facts about differentials on algebraic curves, see e.g.~\cite{Gr}.
Let $X$ be a non-singular projective curve
over $\C$ of genus $g$.
 A differential of the first kind on $X$
is a holomorphic $1$-form. A differential of the second kind is a
meromorphic $1$-form on $X$ whose residues all vanish. A
differential of the third kind on $X$ is a meromorphic $1$-form on
$X$ whose poles are all of first order with residues in $\Z$.
In Section \ref{sect:3} and the subsequent sections
we will relax the condition on the integrality of the residues.
Let
$\psi$ be a differential of the third kind on $X$ that has poles at
the points $P_j$, with residues $c_j$, and is holomorphic elsewhere.
Then the {\it residue divisor} of $\psi$ is
\[
\res(\psi):=\sum_j c_j P_j.
\]
By the residue theorem, the restriction of this divisor to any
component of $X$ has degree $0$.

Conversely, if $D=\sum_j c_j P_j$ is any divisor on $X$ whose
restriction to any component of $X$ has degree $0$, then the Riemann-Roch
theorem and Serre duality imply that there is a differential $\psi_D$ of the third
kind with residue divisor $D$ (see e.g.~\cite{Gr}, p.~233).
Moreover, $\psi_D$ is determined by
this condition up to addition of a differential of the first kind.
Let $U=X\setminus\{ P_j\}$. The canonical homomorphism $H_1(U,\Z)\to
H_1(X,\Z)$ is surjective and its kernel is spanned by the classes of
small circles $\delta_j$ around the points $P_j$. In particular, we
have $\int_{\delta_j} \psi_D =2\pi i c_j$.

Using the Riemann period relations, it can be shown that there is a
unique  differential  of the third kind $\eta_D$ on $X$ with residue
divisor $D$ such that
\[
\Re \left(\int_{\gamma} \eta_D\right) =0
\]
for all $\gamma\in H_1(U,\Z)$. It is called the \emph{canonical
differential of the third kind} associated with $D$. For instance,
if $f$ is a meromorphic function on $X$, then
$df/f$ is a canonical differential of the third kind on $X$ with
residue divisor $\dv(f)$.

Let $\bar\Q\subset\C$ be the algebraic closure of $\Q$ in $\C$.
For the rest of this section we assume that $X$ is defined over $\Q$.
Moreover, we assume that the divisor $D$ is defined over a
number field $F\subset \bar\Q$. Results by Waldschmidt on the transcendence of periods of
differentials of the third kind (see Section 5.2 of \cite{W}, and
Theorem 2 of \cite{Sch}) imply the following theorem.

\begin{theorem}[Scholl]
\label{scholl} If some non-zero multiple of $D$ is a principal
divisor, then $\eta_D$ is defined over $F$. Otherwise, $\eta_D$ is
not defined over $\bar \Q$.
\end{theorem}

We let $J=J(X)$ be the Jacobian of $X$ over $\Q$.
If $k\subset\C$ is a subfield, we denote by $J(k)$ the group of $k$-valued points of $J$. It can be described as the quotient of the group $\Div(X,k)^0$ of degree $0$ divisors on $X$ which are rational over $k$ modulo the subgroup of principal divisors $\dv(f)$ for $f\in k(X)^\times$.
If $k$ is a number field, $J(k)$ is a finitely generated abelian group.

Let $\Omega(X,k)$ denote the
space of differentials of the first kind on $X$ defined over $k$.
We write $\calD(X,k)$ for the group of differentials of the third kind  on $X$ defined over $k$.
We write $\calP(X,k)$ for the subgroup of differentials of the third kind of the form
$df/f$ where $f\in k(X)^\times$ is a rational function defined over $k$.
Moreover, we put
\[
\ccl(X,k)= \calD(X,k)/\calP(X,k).
\]
Associating to a differential of the third kind its residue divisor induces the exact sequence of abelian groups:
\begin{align}
\label{eq:cl}
\xymatrix{ 0\ar[r]& \Omega(X,k) \ar[r]& \ccl(X,k)\ar[r]&  J(k)\ar[r] & 0 }.
\end{align}

Let $K\subset \C$ be a subfield,
and denote by $k\cdot K$ the compositum of $k$ and $K$.
We put $J(k)_K=J(k)\otimes_\Z K$ and $\Div(X,k)^0_K=\Div(X,k)^0\otimes_\Z K$.
We write $\calD(X,k)_K$ for the group of meromorphic differentials on $X$ defined over
$k\cdot K$
whose poles are all of first order
and whose residue divisor belongs to $\Div(X,k)^0_K$.
We write $\calP(X,k)_K$ for the subgroup of differentials which are finite
$K$-linear combinations of differentials of the form
$df/f$ with $f\in k(X)^\times$.
Moreover, we put
\[
\ccl(X,k)_K= \calD(X,k)_K/\calP(X,k)_K.
\]
%
We also have the exact sequence of $K$-vector spaces
\begin{align}
\label{eq:cl2}
\xymatrix{ 0\ar[r]& \Omega(X,k\cdot K) \ar[r]& \ccl(X,k)_K\ar[r]&  J(k)_K\ar[r] & 0 }.
\end{align}

\subsection{Differentials of the third kind on modular curves}
\label{sect:3}

Here we study the action of the Hecke algebra on differentials of the third kind on modular curves.

Let $N$ be a positive integer, and let $X=X_0(N)$ be the modular curve
of level $N$ associated with the group $\Gamma_0(N)$. We write $J$ for
the Jacobian of $X$.  Let $S_2(N)$ be the space of cusp forms of
weight $2$ with respect to $\Gamma_0(N)$.  We identify $S_2(N)$ with
$\Omega(X,\C)$ by the map $G\mapsto \omega_G = 2\pi i \cdot G(z)
\,dz$.
A meromorphic differential $\psi$ on $X$ has a Fourier expansion of the form
\[
\psi= 2\pi i\sum_{n=n_0}^\infty c(n) e^{2\pi i n z} dz.
\]
We refer to $c(n)$ as the $n$-th Fourier coefficient of $\psi$.

Let $k\subset \C$ be a number field. By the $q$-expansion principle, the elements of $\Omega(X,k)$ can be identified with those differentials whose Fourier coefficients are contained in $k$.
The abstract Hecke algebra $\T$ of $\Gamma_0(N)$ acts on $X$ by correspondences, which are defined over $\Q$.
This induces compatible actions of $\T$ on
$\Omega(X,k)$, $\ccl(X,k)$ and $J(k)$.

Throughout this section, let $G\in S_2(N)$ be a (normalized) newform.
We denote by $K=K_G$ the
field of definition of $G$, that is, the totally real number field
generated by the Hecke eigenvalues of $G$.

\begin{proposition}
\label{prop:hecke}
There is an element of\/ $\T\otimes_\Z K$ such that the corresponding Hecke operator on $S_2(N)$ is the orthogonal projection $\C G$.
\end{proposition}

\begin{proof}
Let $T\in \T$. Since $G$ is a newform, we have $T G= \lambda_T G$ with some eigenvalue $\lambda_T\in K$.

We first show that the orthogonal projection to the eigenspace
\[
E(T,\lambda_T)=\{ F\in S_2(N);\; TF=\lambda_T F\}
\]
is given by an element
of $\T\otimes_\Z K$.
To this end, let $P_T(X)\in \Q[X]$ be the characteristic polynomial of the endomorphism of $S_2(N)$ corresponding to $T$.
We write
\[
P_T(X)= (X-\lambda_T)^\mu\cdot Q(X),
\]
where $\mu$ is the algebraic multiplicity of the eigenvalue $\lambda_T$ and $Q\in K[X]$ is a polynomial with $Q(\lambda_T)\neq 0$.
Then $Q(T)\in \T\otimes_\Z K $ acts by multiplication with $Q(\lambda_T)$ on $E(T,\lambda_T)$ and takes any eigenform in the orthogonal complement to zero.
Consequently,
\begin{align}
\label{eq:1}
A(T):=\frac{Q(T)}{Q(\lambda_T)}\in \T\otimes_\Z K
\end{align}
induces the orthogonal projection to $E(T,\lambda_T)$.

If $p$ is a prime, we write $\lambda_p$ for the eigenvalue of $T_p$ corresponding $G$.
Since $G$ is a newform, multiplicity one implies that there exist primes $p_1,\dots, p_r$, such that the common eigenspace of the $T_{p_i}$ is given by
\[
E(T_{p_1},\lambda_{p_1})\cap\dots \cap E(T_{p_r},\lambda_{p_r})=\C G.
\]
The product $A(T_{p_1})\cdots A(T_{p_r})$ of the elements of $\T\otimes_\Z K$  corresponding to $T_{p_i}$ as in \eqref{eq:1} induces the
orthogonal projection to $\C G$.
\end{proof}

Proposition \ref{prop:hecke} implies that the $G$-isotypical
component $J(k)_\C^G$ of $J(k)_\C$ corresponding to $G$
is defined over $K$, that is, it
has a basis consisting of elements of
$J(k)_{K}$.
Analogously, the $G$-isotypical component $\ccl(X,k)_{\C}^G$ of $\ccl(X,k)_{\C}$ is defined over $K$.


If $\alpha\in \calD(X,k)_\C$ and $\beta\in\Omega(X,\C)$, we denote their
$L^2$-pairing by
\[
\langle \alpha,\beta\rangle = \int_X \alpha\wedge \bar \beta.
\]
It is a consequence of Stokes' theorem that the pairing vanishes if $\alpha \in \calP(X,k)_\C$.
Therefore it induces
a pairing of $\ccl(X,k)_\C$ and $\Omega(X,\C)$.  Moreover, the usual
argument shows that the Hecke operators are self adjoint with respect
to this pairing.


\begin{proposition}
\label{prop:1}
Let $D\in \Div(X,k)^0_{K}$ be a
divisor
whose class in $J(k)_{K}$ lies in the $G$-isotypical component.
\begin{itemize}
\item[(i)]
There is a $\zeta_D\in \calD(X,k)_{K}$
with $\res(\zeta_D) = D$ whose class belongs to the $G$-isotypical component of $\ccl(X,k)_{K}$.
\item[(ii)]
The following sequence is exact:
\begin{align*}
\label{eq:clg}
\xymatrix{ 0\ar[r]& (k\cdot K) \omega_G  \ar[r]& \ccl(X,k)_{K}^G\ar[r]&  J(k)_K^G\ar[r] & 0 }.
\end{align*}
\end{itemize}
\end{proposition}

\begin{proof}
i) By the Riemann-Roch
theorem and Serre duality, there is a differential of the third
kind  $\psi_D\in \calD(X,k)_{K}$  with residue divisor $D$.

According to Proposition \ref{prop:hecke}, there is an element $T\in \T\otimes_\Z K$ such that the corresponding Hecke operator on $S_2(N)$ is the orthogonal projection to the $G$-isotypical component.
We have $ \res(T\psi_D) = T D$. Since the class of $D$ is $G$-isotypical, we find that
\[
\res( T\psi_D -\psi_D) =0\in J(k)_K.
\]
In view of \eqref{eq:cl2} there exists an $\alpha\in \calP(X,k)_K$ such that
\[
\beta:= T\psi_D -\psi_D-\alpha\in \Omega(X,k\cdot K).
\]
Observe that $T\beta =0$. In fact, we have
\begin{align*}
\langle T\beta, T\beta\rangle &= \langle \beta, T\beta\rangle\\
&=\langle T\psi_D , T\beta\rangle -\langle\psi_D, T\beta\rangle -\langle\alpha, T\beta\rangle \\
&= \langle \psi_D , T^2\beta\rangle -\langle\psi_D, T\beta\rangle \\
&=0.
\end{align*}

We claim that
\[
\zeta_D:=\psi_D+\beta
\]
has the required properties. To see this we note that
$\res(\zeta_D)= \res(\psi_D)=D$, and
\begin{align*}
T\zeta_D-\zeta_D&= T\psi_D+T\beta-\psi_D-\beta\\
&=T\psi_D-\psi_D-\beta\\
&=\alpha.
\end{align*}
Consequently, the class of $\zeta_D$ belongs to the $G$-isotypical component of $\ccl(X,k)_K$.

ii) The second statement is a direct consequence of the first.
\end{proof}

\begin{proposition}
\label{prop:3}
Let $D\in \Div(X,k)^0_K$ be a divisor  whose class in $J(k)_{K}$ lies in the $G$-isotypical component.
Then there is a unique $\zeta_D\in \calD(X,k)_{K}$ with the following properties:
\begin{itemize}
\item[(i)]
The residue divisor is given by $\res(\zeta_D) = D$.
\item[(ii)] The class of $\zeta_D$ belongs to the $G$-isotypical
component of $\ccl(X,k)_{K}$.
\item[(iii)]
The first Fourier coefficient of $\zeta_D$ vanishes.
\end{itemize}
\end{proposition}

\begin{proof}
According to Proposition \ref{prop:1}, there exists a differential $\psi_D\in \calD(X,k)_{K}$ satisfying (i) and (ii).
Then for any $a\in k\cdot K$, the linear combination
\[
\psi_D+a  \omega_G
\]
satisfies (i) and (ii) as well. We write $c(1)$ for the first Fourier
coefficient of $\psi_D$.  Since the first Fourier coefficient of
$\omega_G$ is equal to $1$,  we obtain a
differential $\zeta_D$ with the required properties by choosing
$a=-c(1)$.

If $\tilde \zeta_D$ is another differential satisfying (i)--(iii), then $\delta:= \tilde \zeta_D-\zeta_D$ has residue divisor $0$
and is therefore holomorphic.  Because of property (ii) and
multiplicity one for the newform $G$, the differential $\delta$ is a
multiple of $\omega_G$.  Since the first coefficient of $\delta$
vanishes, we find that $\delta=0$, and thereby $ \tilde \zeta_D=\zeta_D$.
\end{proof}

\begin{definition}
We call the differential $\zeta_D$ in Proposition \ref{prop:3} the \emph{normalized differential of the third kind} associated with $D$.
\end{definition}

\subsection{Comparing canonical and normalized differentials of the third kind}

Let $H_1(X,\R)=H_1(X,\Z)\otimes_\Z \R$ be the first homology of $X$ with real coefficients. The map $\H\to \H$, $z\mapsto -\bar z$ induces a real analytic automorphism $\sigma$ of the modular curve $X$. It induces an $\R$-linear involution on $H_1(X,\R)$. We let $H_1^\pm(X,\R)$ be the eigenspaces corresponding to the eigenvalues $\pm 1$. The map $\sigma$ also induces an $\R$-linear involution on $\calD(X,k)_\C$, which we also denote by $\sigma$. If we view the elements of $\calD(X,k)_\C$ as meromorphic modular forms of weight $2$, then $\sigma$ corresponds to the involution given by complex conjugation of the Fourier coefficients.
We denote the eigenspaces associated with the eigenvalues $\pm 1$ by  $\calD(X,k)_\C^\pm$. Analogously, we write $\Omega(X,\C)^\pm$ for the eigenspaces of $\sigma$ on the holomorphic differentials.
Note that
$\Omega(X,\C)^\pm$ can be identified with the subspace of elements of $S_2(N)$ with real (respectively imaginary) Fourier coefficients.

We consider the bilinear pairing
\begin{align}
\label{eq:pair0}
\Omega(X,\C) \times H_1(X,\R)\longrightarrow \C,
\quad (\omega,C)\mapsto \langle\omega,C\rangle = \int_C\omega.
\end{align}
It is well known that the restriction
\begin{align}
\label{eq:pair1}
\Omega(X,\C)^\pm \times H_1^\pm(X,\R)\longrightarrow \R
\end{align}
takes real values and defines a non-degenerate $\R$-bilinear pairing (see e.g. \cite[Chapter~2]{Cr}).
The Hecke algebra $\T$ acts on $H_1^\pm(X,\R)$, and for $T\in\T$ we
have $\langle T\omega ,C\rangle= \langle \omega, TC\rangle$.  If $G\in
S_2(N)$ is a newform, then we write $H_1^{\pm,G}(X,\R)$ for the
corresponding $G$-isotypical component of $H_1^\pm(X,\R)$. It is a
one-dimensional subspace.  We fix a generator $C_G^\pm\in
H_1^{\pm,G}(X,\R)$.  Proposition \ref{prop:hecke} implies that
$C_G^\pm$ can be chosen in $H_1^\pm(X,K)$, but we do not require
this at the moment.

We extend the pairing \eqref{eq:pair1} to an $\R$-bilinear pairing
\begin{align}
\label{eq:pair2}
\calD(X,k)_\R\times H_1(X,\R)\longrightarrow \R, \quad (\psi,C)\mapsto
\Re\langle\psi,C\rangle = \Re \int_C\psi.
\end{align}
Here, to compute the integral, we take a representative for $C$ which
is disjoint from the support of $\res(\psi)$. By the residue theorem,
the value of the pairing is independent of the choice of such a
representative. By Stokes' theorem, the pairing vanishes if $\psi\in
\calP(X,k)_\R$.

If $D\in \Div(X,k)^0$, then the canonical differential of the third kind corresponding to $D$ is the unique differential $\eta_D\in \calD(X,k)_\R$  such that $\res(\eta_D)=D$ and $\Re\langle \eta_D, C\rangle = 0$ for all $C\in H_1(X,\R)$. The following proposition compares canonical and normalized differentials.
Let $\Div(X,k)^{0,\pm}_K$ be the eigenspace of $\sigma$ on $\Div(X,k)^{0}_K$ corresponding to the eigenvalue $\pm 1$.

\begin{theorem}
\label{thm:comp0}
Let $D\in \Div(X,k)^{0,\pm}_K$ be a divisor whose class in $J(k)_{K}$ lies in the $G$-isotypical component. Let $\eta_D$ be the canonical differential of the third kind associated with $D$, and let $\zeta_D$ be the normalized differential of the third kind associated with $D$.
Then we have
\[
\eta_D= \zeta_D -\frac{\Re\langle \zeta_D,C_G^\pm \rangle}{\langle \omega_G,C_G^\pm\rangle}\cdot \omega_G.
\]
\end{theorem}

\begin{proof}
Since $\res(\eta_D)=\res(\zeta_D)$, the difference $\delta:=\eta_D-\zeta_D$ is a holomorphic differential in  $\Omega(X,k\cdot\R)$.
It is easily checked that $\sigma(\eta_D)=\eta_{\sigma(D)}$ and $\sigma(\zeta_D)=\zeta_{\sigma(D)}$. Hence the assumption on $D$ implies that
$\delta\in \Omega(X,k\cdot\R)^\pm$.
To determine $\delta$,  we compute its pairing with  $H_1^\pm(X,\R)$.

According to Proposition \ref{prop:hecke}, there is an element $T\in \T\otimes_\Z K_G$ such that the corresponding Hecke operator on $S_2(N)$ is the orthogonal projection to the $G$-isotypical component.
For  $C\in H_1^\pm(X,\R)$ we have
\begin{align*}
\langle T\delta ,C\rangle &= \Re\langle \eta_D ,T C\rangle -\Re\langle T\zeta_D ,C\rangle\\
&= -\Re\langle \zeta_D ,C\rangle\\
&= \langle \delta ,C\rangle .
\end{align*}
Here we have used that $\eta_D$ is a canonical differential and that $\zeta_D$ is $G$-isotypical. Since the pairing \eqref{eq:pair1} is non-degenerate, we find that $T\delta = \delta$. So $\delta = a\omega_G$ for some $a\in \C$. Using the action of $\sigma$ we see that $a\in \R$ if $D\in \Div(X,k)^{0,+}_K$, and  $a\in i\R$ if $D\in \Div(X,k)^{0,-}_K$.

To determine $a$, we compute the pairing with $C_G^\pm$. We obtain
\begin{align*}
a\langle \omega_G ,C_G^\pm\rangle &= \langle \delta ,C_G^\pm\rangle = -\Re\langle \zeta_D ,C_G^\pm\rangle.
\end{align*}
Since  the pairing \eqref{eq:pair1} is non-degenerate, the quantity $\langle \omega_G ,C_G^\pm\rangle$ does not vanish, and therefore
\[
\delta=-\frac{\Re\langle \zeta_D,C_G^\pm\rangle}{\langle \omega_G,C_G^\pm\rangle}\cdot \omega_G.
\]
This concludes the proof of the proposition.
\end{proof}

The following corollary gives an interpretation of  the first Fourier coefficient of $\eta_D$ as the quotient of two periods of algebraic differentials on $X$.

\begin{corollary}
\label{cor:comp}
Let $D\in \Div(X,k)^{0,\pm}_K$ be a divisor whose class in $J(k)_{K}$ lies in the $G$-isotypical component. The first Fourier coefficient of $\eta_D$ is given by
\begin{align}
\label{eq:quant}
-\frac{\Re\langle \zeta_D,C_G^\pm\rangle}{\langle \omega_G,C_G^\pm\rangle}.
\end{align}
\end{corollary}

\begin{remark}
  When $K=\Q$, then Theorem \ref{scholl} implies that the quantity \eqref{eq:quant} is
  algebraic if and only if the image of $D$ in $J(k)_{K}$
  vanishes.
\end{remark}


The next corollary
  gives an alternative characterization of the normalized differential
  of the third kind corresponding to $D$.

\begin{corollary}
\label{cor:altchar}
Let $D\in \Div(X,k)^{0}_K$ be a degree $0$ divisor  whose class in $J(k)_{K}$ lies in the $G$-isotypical component.
Then there is a unique $\zeta_D\in \calD(X,k)_{\R}$ with the following properties:
\begin{itemize}
\item[(i)]
The residue divisor is given by $\res(\zeta_D) = D$.
\item[(ii)] We have $\Re\langle \zeta_D,C\rangle=0$ for all $C\in H_1(X,\R)$ with $\langle \omega_G, C\rangle=0$.
\item[(iii)]
The first Fourier coefficient of $\zeta_D$ vanishes.
\end{itemize}
\end{corollary}

In view of Proposition \ref{prop:3}, such a $\zeta_D$ will automatically be contained in $\calD(X,k)_{K}$.

\section{Harmonic Maass forms and Heegner divisors}

Here we relate periods of normalized differentials of the third kind associated with Heegner divisors to coefficients of harmonic Maass forms.
We use $\tau=u+iv$ with $u,v\in \R$ as a standard variable in the upper complex half plane $\H$.
Let $N$ be a positive integer and let $X=X_0(N)$ be the (projective) modular curve of level $N$ associated with $\Gamma_0(N)$. As before, we write $J$ for the Jacobian of $X$.

\subsection{Twisted Heegner divisors}

Let $d,h\in \Z$ be integers such that $d\equiv h^2\pmod{4N}$.
Let $\calQ_{d,h}$ be the set of integral binary quadratic forms $Q=[Na,b,c]$ of discriminant $b^2-4Nac=d$ for which $b\equiv h\pmod{2N}$. Here $a,b,c\in \Z$.
The group $\Gamma_0(N)$ acts on $\calQ_{d,h}$ with finitely many orbits.
If $d<0$, we write $\alpha_Q$ for the Heegner point associated with $Q$, that is, the zero in $\H$ of the polynomial $aX^2+bX+c\in \Z[X]$.
In this case we write $w_Q$ for the order of the (finite) stabilizer of $Q$ in $\Gamma_0(N)$.

Let  $\Delta\in \Z$ be a {\em fundamental} discriminant and $r$ be an integer such that $\Delta \equiv r^2\pmod{4N}$.
If $\Delta\mid d$ and $d/\Delta$ is also a square modulo $4N$, then we have a `genus character' $\chi_\Delta$ on $\calQ_{d,h}$ given by
\[
\chi_{\Delta}(Q)=\chi_{\Delta}([Na,b,c]):=
\begin{cases}
\leg{\Delta}{n},&\text{if $\operatorname{gcd}(a,b,c,\Delta)=1$,}\\
0,& \text{otherwise}.
\end{cases}
\]
Here $n$ is any integer prime to $\Delta$ represented
by one of the quadratic forms $[N_1a,b,N_2 c]$ with $N_1N_2=N$ and
$N_1,N_2>0$. Such an integer exists and the value of $\leg{\Delta}{n}$ is independent of the choice of $N_1$, $N_2$, and $n$
(see Section 1.2 of \cite{GKZ}, and also Section 1 of \cite{Sk2}).

Now assume that $d$ and $\Delta$ are discriminants with opposite sign.
Then $\Delta d$ is a negative discriminant which is a square modulo $4N$.
We define the
{\it twisted Heegner divisor} $Z_{\Delta,r}(d,h)$ by
\begin{equation}
Z_{\Delta,r}(d,h) := \sum_{Q\in \calQ_{\Delta d,rh}/\Gamma_0(N)}
\chi_\Delta (Q)\cdot  \frac{\alpha_Q}{w_Q}.
\end{equation}
Then $Z_{\Delta,r}(d,h)$ is a divisor on $X$ defined over
$\Q(\sqrt{\Delta})$ (see \cite{BO}, Lemma 5.1). We have that $\sigma Z_{\Delta,r}(d,h)= \sgn(\Delta)Z_{\Delta,r}(d,h)$. The divisor $Z_{\Delta,r}(d,h)$ has degree $0$ unless $\Delta=1$, in which case the degree is the Hurwitz class number $H(d)$.
The divisor
\begin{equation}
y_{\Delta,r}(d,h):=Z_{\Delta,r}(d,h) -
\deg(Z_{\Delta,r}(d,h))\cdot \infty
\end{equation}
has degree $0$ for any $\Delta$. It defines a point in $J(\Q(\sqrt{\Delta}))$.

\subsection{Harmonic Maass forms}

We write $\Mp_2(\R)$ for the metaplectic two-fold cover of
$\Sl_2(\R)$, realized as the group of pairs $(M,\phi(\tau))$,
where $M=\kabcd\in\Sl_2(\R)$ and $\phi:\H\to \C$ is a holomorphic
function with $\phi(\tau)^2=c\tau+d$, see e.g. \cite{Bo1}, \cite{Br}.
We denote  the inverse image of
$\Gamma:=\Sl_2(\Z)$ under the covering map by
$\tilde\Gamma:=\Mp_2(\Z)$. It is well known that $\tilde\Gamma$ is
generated by $T:= \left( \kzxz{1}{1}{0}{1}, 1\right)$, and $S:=
\left( \kzxz{0}{-1}{1}{0}, \sqrt{\tau}\right)$.

Let $N$ be a positive integer. There is a `Weil representation' $\rho$ of $\tilde \Gamma$ on
$\C[Z/2N\Z]$, the group ring of the cyclic group of order $2N$.
For a coset $h\in \Z/2N\Z$ we denote by $\frake_h$ the corresponding standard basis vector of $\C[\Z/2N\Z]$.
In terms of the generators $T$ and $S$ of $\tilde \Gamma$, the representation $\rho$ is given by
\begin{align}
\label{eq:weilt}
\rho(T)(\frake_h)&=e\left(\frac{h^2}{4N}\right)\frake_h,\\
\label{eq:weils}
\rho(S)(\frake_h)&=
\frac{1}{\sqrt{2iN}} \sum_{h' \; (2N)} e\left(-\frac{hh'}{2N}\right)
 \frake_{h'}.
\end{align}
Here the sum runs through the elements of $\Z/2N\Z$ and we have put $e(a)=e^{2\pi i a}$.
Note that $\rho$ is the Weil representation associated with the one-dimensional positive definite lattice $K=(\Z,Nx^2)$ in the sense of \cite{Bo1}, \cite{Br}, \cite{BO}.
It is unitary with respect to the standard scalar product.

If $k\in \frac{1}{2}\Z$, we write $M^!_{k,\rho}$ for the space of $\C[\Z/2N\Z]$-valued weakly holomorphic modular forms of weight $k$ for $\tilde \Gamma$ with representation $\rho$. The subspace of holomorphic modular forms (respectively cusp forms) is denoted by $M_{k,\rho}$ (respectively $S_{k,\rho}$).
According to \cite[Chapter 5]{EZ}, the space $M_{k,\bar\rho}$ is isomorphic to $J_{k+1/2,N}$, the space of holomorphic Jacobi forms of weight $k+1/2$ and index $N$. According to \cite{Sk1} and \cite{SZ}, $M_{k,\rho}$ is isomorphic to $J_{k+1/2,N}^{skew}$, the space of skew holomorphic Jacobi forms of weight $k+1/2$ and index $N$.

Let $k\in \frac{1}{2}\Z$.
A twice continuously differentiable function $f:\H\to
\C[\Z/2N\Z]$ is called a {\em harmonic Maass form} (of weight $k$ with
respect to $\tilde \Gamma$ and $\rho$) if it satisfies:
\begin{enumerate}
\item[(i)]
$f(M\tau) = \phi(\tau)^{2k}\rho(M,\phi) f(\tau)$
for all $(M,\phi)\in \tilde\Gamma$;
\item[(ii)]
$\Delta_k f=0$,
\item[(iii)]
there is a $\C[\Z/2N\Z]$-valued Fourier polynomial
\[
P_f(\tau)=\sum_{h\;(2N)}\sum_{n\in \Z_{\leq 0}} c^+(n,h) q^{\frac{n}{4N}} \frake_h
\]
such that $f(\tau)-P_f(\tau)=O(e^{-\eps v})$ as $v\to \infty$ for
some $\eps>0$.
\end{enumerate}
%
Here we have that
\begin{align*}
\Delta_k := -v^2\left( \frac{\partial^2}{\partial u^2}+
\frac{\partial^2}{\partial v^2}\right) + ikv\left(
\frac{\partial}{\partial u}+i \frac{\partial}{\partial v}\right)
\end{align*}
is the usual weight $k$ hyperbolic Laplace operator (see \cite{BF}).
The Fourier polynomial $P_f$  is called the {\em principal part} of
$f$. We denote the vector space of these harmonic Maass forms
by  $H_{k,\rho}$ (it was called $H^+_{k,\rho}$ in \cite{BF}).
Any $f\in H_{k,\rho}$ has a Fourier expansion of the form
\begin{align}
\label{eq:fourier}
f(\tau)= \sum_{h\;(2N)}\sum_{\substack{n\in \Z\\ n\gg-\infty}} c^+(n,h) q^{\frac {n}{4N}}\frake_h
+ \sum_{h\; (2N)}\sum_{\substack{n\in \Z\\
n< 0}} c^-(n,h) \Gamma(1-k,4\pi|n|v) q^{\frac{n}{4N}} \frake_h.
\end{align}

Recall that there is an antilinear  differential operator $\xi=
\xi_k:H_{k,\rho}\to S_{2-k,\bar\rho}$, defined by
\begin{equation}
\label{defxi} f(\tau)\mapsto \xi(f)(\tau):=2iv^k\overline{\frac{\partial f}{\partial \bar \tau}}.
\end{equation}
Here $\bar\rho$ denotes the dual of the representation $\rho$ (which can be identified with the complex conjugate with respect to the standard basis).
The map $\xi$ is surjective and its kernel is $M^!_{k,\rho}$.

Let  $\Delta\in \Z$ be a {\em fundamental} discriminant and $r$ be an integer such that $\Delta \equiv r^2\pmod{4N}$.
If $\Delta>0$ we let $\rho'=\rho$, and if $\Delta<0$ we let $\rho'=\bar\rho$.
Let $f\in H_{1/2,\rho'}$ be a harmonic Maass form and denote its Fourier coefficients by $c^\pm(n,h)$ as in \eqref{eq:fourier}.
Assume that the coefficients $c^+(n,h)$ are real.
We associate a degree $0$ divisor with $f$ by putting
\[
y_{\Delta,r}(f)=\sum_{h\;(2N)}\sum_{n\in \Z_{\leq 0}} c^+(n,h) y_{\Delta,r}(n,h).
\]
Let  $\eta_{\Delta,r}(f)$ be the canonical differential of the third kind on $X$ associated with $y_{\Delta,r}(f)$. In \cite{BO} it is constructed as a regularized theta lift of $f$ and its Fourier expansion is computed. According to \cite[Theorem 5.4]{BO} we have:

\begin{theorem}
\label{cdfourier}
Assume that $\Delta\neq 1$.
The canonical differential  associated with $y_{\Delta,r }(f)$ is given by
\begin{align*}
\eta_{\Delta,r }(z,f)
&= -
\sgn(\Delta)\sqrt{\Delta}\sum_{\substack{n\geq 1}}\sum_{d\mid n}
\frac{n}{d}\leg{\Delta}{d}
c^+(\tfrac{|\Delta|n^2}{d^2},\tfrac{rn}{d}) e(nz)\cdot
2\pi i\,dz.
\end{align*}
\end{theorem}

Note that for $\Delta=1$ the same statement is true if $y_{\Delta,r}(f)$ is corrected by some degree $0$ divisor supported at
the cusps and the Fourier expansion on the right hand side is corrected by an additional constant term (given by the Weyl vector
associated with $f$ at the cusp $\infty$).

\subsection{Fourier coefficients of harmonic Maass forms}
\label{sect:4.3}

We now recall one of the main results of \cite{BO} on the algebraicity of the coefficients of harmonic Maass forms.
Let $G\in
S_{2}(N)$ be a normalized newform (in particular a
common eigenform of all Hecke operators) of weight $2$ and write
$K=K_G$ for the field of definition of $G$.
Let $\eps=\eps_G\in \{\pm 1\}$ be the eigenvalue of the Fricke involution on $G$, that is, $G\mid W_N = \eps G$. Then the sign
of the functional equation of the Hecke $L$-function of $G$ is equal to $-\eps$.
%

If
$\eps=+1$ we put $\rho'=\rho$,
and
if $\eps=-1$ we put $\rho'=\bar \rho$.
There is a newform
$g\in S_{3/2,\bar\rho'}$ mapping to $G$ under the Shimura
correspondence. It is well
known that we may normalize $g$ such that all its coefficients
are contained in $K$. According to \cite[Lemma 7.3]{BO}, there is a
harmonic  Maass form $f\in H_{1/2, \rho'}$ whose principal part has coefficients in $K$ with the property that
\[
\xi_{1/2}(f)=\|g\|^{-2} g.
\]
Possibly replacing $f$ by $\lambda f$ and $g$ by $g/\lambda$ for a suitable positive integer $\lambda$, we may
actually assume that the principal part of $f$ has coefficients in the ring $\calO_K$ of integers of $K$.
This form is then unique up to addition of a weakly holomorphic form in
$M^!_{1/2, \rho'}$ whose principal part has coefficients in $\calO_K$.
According to \cite[Theorem 7.5]{BO}, the divisor $y_{\Delta,r}(f)$ defines a point in the $G$-isotypical component $J(\Q(\sqrt{\Delta}))_K^G$.

\begin{theorem}[see \cite{BO}]
\label{Lval2}
Let $G\in S_{2}(N)$ be a normalized newform with $G\mid W_N = \eps G$.
Let $g\in S_{3/2,\bar\rho'}$ and $f\in H_{1/2, \rho'}$ be as above. Denote the Fourier coefficients of $f$ by $c^\pm(n,h)$ for $n\in \Z$ and $h\in \Z/2N\Z$. Then the following are true:
\begin{enumerate}
\item
If $\Delta\neq 1$ is a fundamental
discriminant and $r\in \Z$ such that $\Delta\equiv r^2\pmod{4N}$ and $\varepsilon_G \Delta <0$, then
$$
L(G,\chi_{\Delta},1)=8\pi^2\|G\|^2 \|g\|^2 \sqrt{\frac{|\Delta|}{N}}\cdot
c^{-}(\varepsilon_G\Delta,r)^2.
$$
\item
If $\Delta\neq 1$ is a fundamental
discriminant and $r\in \Z$ such that $\Delta\equiv r^2\pmod{4N}$ and $\varepsilon_G \Delta >0$, then
$$
L'(G,\chi_{\Delta},1)=0 \quad \Longleftrightarrow\quad c^{+}(\varepsilon_G\Delta,r)\in \bar \Q\quad   \Longleftrightarrow\quad c^{+}(\varepsilon_G\Delta,r)\in K.$$
\end{enumerate}
\end{theorem}

When $S_{1/2,\rho'}=\{0\}$ the above result also holds for $\Delta=1$, see also \cite[Remark 18]{BO}. This is for instance the case when $N$ is a prime.
If $N$ is a prime and $\eps=1$, then the space $H_{1/2, \rho'}$ can be identified with a space of scalar valued modular forms satisfying a Kohnen plus space condition.

We now combine the results of Section \ref{sect:3}, in particular Corollary~\ref{cor:comp}, with Theorem~\ref{cdfourier} to obtain an exact formula form the coefficients $c^{+}(\varepsilon_G\Delta,r)$ of the holomorphic part of $f$.

\begin{theorem}
\label{thm:comp}
Let $G\in S_{2}(N)$ be a normalized newform, and let $\eps\in \{\pm 1\}$ such that  $G\mid W_N = \eps G$.
Let $g\in S_{3/2,\bar\rho'}$ and $f\in H_{1/2, \rho'}$ be as above.
Assume that $\Delta\neq 1$ is a fundamental
discriminant and $r\in \Z$ such that $\Delta\equiv r^2\pmod{4N}$ and $\eps \Delta >0$, and let $\zeta_{\Delta,r}(f)\in \calD(X,\Q(\sqrt{\Delta}))_K^{\eps}$ be the normalized differential of the third kind associated with $y_{\Delta,r}(f)$. Then
\[
c^{+}(\eps \Delta,r)=\eps\cdot \frac{\Re\langle \zeta_{\Delta,r}(f),\,C_G^{\eps}\rangle}{ \sqrt{\Delta}\langle \omega_G,\,C_G^{\eps}\rangle}.
\]
\end{theorem}

\begin{proof}
According to Theorem \ref{cdfourier}, the first Fourier coefficient of $\eta_{\Delta,r}(f)$ is given by
\[
-\sgn(\Delta)\sqrt{\Delta} c^+(|\Delta|,r).
\]
On the other hand, Corollary \ref{cor:comp} says that this coefficient is equal to
\[
-\frac{\Re\langle \zeta_{\Delta,r}(f),\,C_G^\eps\rangle}{\langle \omega_G,\,C_G^\eps\rangle}.
\]
This implies the assertion.
\end{proof}

\begin{remark}
Note that normalized differentials of the third kind associated with (twisted) Heegner divisors can be constructed explicitly using (a twisted version of) the additive regularized theta lift \cite[Theorem 14.3]{Bo1} from weakly holomorphic modular forms in $M^!_{3/2,\rho'}$ to $\Orth(2,1)$.
\end{remark}

\section{Rational elliptic curves}

In this section we consider the implications of Theorem \ref{thm:comp} for rational elliptic curves.
Let
\[
E: y^2=4x^3-g_2 x-g_3
\]
with $g_2,g_3\in \Q$ be an elliptic curve over $\Q$ of conductor $N$.
Let $X=X_0(N)$ be the modular curve of level $N$ viewed as a projective algebraic curve over $\Q$.
In view of the work of Wiles et al.~\cite{Wi}, \cite{BCDT}, there exists a modular parameterization
\[
\phi: X\longrightarrow E
\]
over $\Q$, that is, a surjective morphism of algebraic curves over $\Q$ that factors through the Jacobian $J$ of $X$.
We may assume that the cusp $\infty$ of $X$ is mapped to the identity element of $E$.

Let
$\omega_E=dx/y$ be the N\'eron differential on $E$.
Then multiplicity one implies  that
\[
\phi^*(\omega_E)=c_E\cdot\omega_G,
\]
where $G\in S_2(N)$ is the normalized newform corresponding to $E$, and $c_E$ denotes the Manin constant.
According to \cite[Proposition 2]{Ed}, $c_E$ is a non-zero integer. Let $\sigma$ be the automorphism of $E$ induced by complex conjugation.
We let $C_E^\pm$ be non-zero elements of $H_1(E,\Z)$ such that $\sigma C_E^\pm = \pm C_E^\pm$. Then $H_1(E,\Q)= \Q C_E^+\oplus \Q C_E^-$. Note that for $C_E^+$ we can take for instance $E(\R)$.
We denote the corresponding periods by
\[
\langle \omega_E,C_E^\pm \rangle = \int_{C_E^\pm }\omega_E.
\]

Let $k$ be a number field. If $Q=(x_0,y_0)\in E(k)$, then
\begin{align}
\label{eq:beta}
\beta(Q):=\frac{y_0}{2(x-x_0)} \cdot \omega_E\in \calD(E,k)_\Q
\end{align}
is a differential of the third kind on $E$ with residue divisor $\frac{1}{2}((Q)-(-Q))$.
Its pullback $\phi^*\beta(Q)$ belongs to $\calD(X,k)_\Q^G$.

Let $\eps=\eps_G$ be the eigenvalue of the Fricke involution on $G$.
Let $f\in H_{1/2, \rho'}$ be a harmonic Maass form corresponding to
$G$ as in Section \ref{sect:4.3}. We may assume that
the principal part of $f$ has integral coefficients. Such an $f$ is
uniquely determined up to addition of a weakly holomorphic form with
integral principal part.  Let $\Delta\neq 1$ be a fundamental
discriminant and $r\in \Z$ such that $\Delta\equiv r^2\pmod{4N}$ and
$\eps \Delta >0$.  Then $y_{\Delta,r}(f)$ defines a
$\Q(\sqrt{\Delta})$-rational point of $J$ which is $G$-isotypical. By
means of the modular parameterization $\phi$, we obtain a point
$P_{\Delta,r}(f)$ in $E(\Q(\sqrt{\Delta}))$.  It is mapped to its
negative under the non-trivial automorphism of $\Q(\sqrt{\Delta})$.
We briefly write $\beta_{\Delta,r}(f)$ for the differential of the
third kind $\beta(P_{\Delta,r}(f))$ corresponding to $P_{\Delta,r}(f)$
as in \eqref{eq:beta}.

\begin{theorem}
\label{thm:ec}
Let the notation be as above.
Then
\[
 c^{+}(\eps\Delta,r) - \frac{\eps c_E }{ \sqrt{\Delta}\langle \omega_E,\,C_E^\eps\rangle }\cdot \Re\int_{C_E^\eps}\beta_{\Delta,r}(f)\in \Q.
\]
\end{theorem}

\begin{proof}
We put $k=\Q(\sqrt{\Delta})$.
According to Theorem \ref{thm:comp} we have
\begin{align}
\label{eq:ec1}
c^{+}(\eps\Delta,r)=\eps\cdot \frac{\Re\langle \zeta_{\Delta,r}(f),\,C_G^\eps\rangle}{ \sqrt{\Delta}\langle \omega_G,\,C_G^\eps\rangle}.
\end{align}
We chose the generator $C_G^\eps\in H_1^{\eps,G}(X,\R)$ such that $\phi_*(C_G^\eps)=C_E^\eps$.
Then we have
\[
\langle \omega_G,\,C_G^\eps\rangle = \frac{1}{c_E}\int_{C_G^\eps}\phi^*(\omega_E)
= \frac{1}{c_E}\langle \omega_E,\,C_E^\eps\rangle .
\]

The pullback $\phi^*\beta_{\Delta,r}(f)$ belongs to
$\calD(X,k)^G_\Q$. The same is true for $\zeta_{\Delta,r}(f)$.
Moreover, by construction, the residue divisors of both differentials
define the same point in the Jacobian.
Hence, according to Proposition \ref{prop:1}, there is an $A_0\in k$ and a $\gamma \in \calP(X,k)_\Q$ such that
\[
\zeta_{\Delta,r}(f)-\phi^*\beta_{\Delta,r}(f)=A_0\omega_G+\gamma.
\]
The non-trivial automorphism of $k$ takes the left hand side to its negative. Hence the same must be true for the right hand side. In particular, we have $A_0=A_1/\sqrt{\Delta} $ for some $A_1\in \Q$.
We obtain  that
\begin{align*}
\Re\langle \zeta_{\Delta,r}(f),\,C_G^\eps\rangle &= \Re\langle\phi^*\beta_{\Delta,r}(f),\,C_G^\eps\rangle+\frac{A_1}{\sqrt{\Delta}}\langle\omega_G,\,C_G^\eps\rangle\\
&=\Re \int_{C_E^\eps}\beta_{\Delta,r}(f)+\frac{A_1\langle \omega_E,\,C_E^\eps\rangle }{c_E\sqrt{\Delta}}.
\end{align*}
Inserting this into \eqref{eq:ec1}, we find that
\[
c^{+}(\eps\Delta,r) =\frac{\eps c_E}{\sqrt{\Delta}\langle \omega_E,\,C_E^\eps\rangle}\cdot
\Re \int_{C_E^\eps}\beta_{\Delta,r}(f)+\frac{\eps A_1}{\Delta} .
\]
This concludes the proof of the theorem.
\end{proof}

The map $\varphi: E\to E_\Delta$, $(x,y)\mapsto (\Delta x, \Delta^{3/2}y)$, defines an isomorphism over $\Q(\sqrt{\Delta})$ of $E$ and its quadratic twist
\[
E_\Delta: v^2=4u^3-\Delta^2g_2 u-\Delta^3 g_3.
\]
We chose $C_E^\eps$ such that $\varphi_*(C_E^\eps)=E_\Delta(\R)$.
The image under $\varphi$ of $P_{\Delta,r}(f)
\in E(\Q(\sqrt{\Delta}))$ defines a rational point $(u_0,v_0)
\in E_\Delta(\Q)$.
If
\[
\alpha_{\Delta,r}(f)=\frac{v_0}{2(u-u_0)}
\cdot \omega_{E_\Delta}\in \calD(E_\Delta,\Q)
\]
denotes the corresponding differential of the third kind as in \eqref{eq:beta}, then
\[
\int_{E_{\Delta}(\R)} \alpha_{\Delta,r}(f) = \int_{C_E^\eps} \beta_{\Delta,r}(f).
\]
Moreover, the real period $\Omega(E_\Delta)= \int_{E_\Delta(\R)}\omega_{E_\Delta}$ of $E_\Delta$ is equal to
$\langle \omega_E, C_E^\eps\rangle /\sqrt{\Delta}$.
We obtain the following corollary.

\begin{corollary}
\label{cor:ec}
Let the notation be as above. Then
\[
\Delta c^{+}(\eps\Delta,r) - \frac{\eps c_E }{\Omega(E_\Delta)}\cdot \Re\int_{E_\Delta(\R)}\alpha_{\Delta,r}(f)\in \Q.
\]
\end{corollary}

We now discuss how this result could probably refined to obtain that the
difference of $\Delta c^+(\varepsilon\Delta,r)$ and a period integral
similar
to the one above is contained in $\frac{1}{4}\Z$.
%
Let
\begin{align}
\label{eq:weierfine}
W: y^2+A_1 xy + A_3 y = x^3+A_2 x^2 + A_4 x +A_6
\end{align}
be a minimal Weierstrass model for $E$ over $\Z$.
Let $\calW\subset \P_\Z^2$ be the closed subscheme defined by $W$. Let $\omega_W=\frac{dx}{2y+A_1x+A_3}$ be the N\'eron differential  on $W$.

\begin{lemma}
\label{lem:wm}
Let $P\in W(\Q)$ and denote by $\calP\in \calW(\Z)$ the Zariski closure of $P$.
There exists a rational section $\underline\alpha(P)$ of the sheaf $\Omega_{\calW/\Z}^1$ of relative differentials such that
\begin{itemize}
\item[(i)]
the induced differential on $\calW_\Q$ has residue divisor $\frac{1}{2}((P)-(-P))$,
\item[(ii)]
the differential $4\underline\alpha(P)$ is regular on $\calW$ up to first order poles on $\calP$ and $-\calP$.
\end{itemize}
The differential $\underline\alpha(P)$ is unique up to addition of an integral multiple of $\frac{1}{4}\omega_W$.
\end{lemma}

\begin{proof}
We write $P=(x_0,y_0)$ with $x_0,y_0\in\Q$.
For any $t\in \Q$, the rational differential
\begin{align}
\label{eq:l0}
\delta_t(P)= \left( \frac{2y_0+A_1 x_0 +A_3}{x-x_0}+t\right)\omega_W
\end{align}
has residue divisor $2(P)-2(-P)$. It defines a rational section of $\Omega_{\calW/\Z}^1$ with first order poles on $\calP$ and $-\calP$ and possible additional poles along vertical divisors of $\calW$. We show that $t$ can be chosen such that there are no vertical poles. Then we can take $\underline\alpha(P)=\frac{1}{4}\delta_t(P)$.

First, notice that  $\frac{A_1 x_0 +A_3}{x-x_0}\omega_W$ has no vertical poles. Therefore it suffices to find a $t\in \Q$ such that
\[
\beta_t(P):= \left(\frac{2y_0}{x-x_0}+t\right)\omega_W
\]
has no vertical poles.
We write $x_0=a/b$ with $a,b\in \Z$ coprime and $y_0=c/d$ with $c,d\in \Z$ coprime.
Since $P$ satisfies \eqref{eq:weierfine}, we have
\begin{align}
\label{eq:l1}
b^3c^2+ A_1 ac b^2 d +A_3 c b^3 d = a^3 d^2 + A_2 a^2 b d^2  + A_4 a b^2 d^2 + A_6 b^3 d^2.
\end{align}
It is easily seen that $b$ divides $d$, so $d=bd'$ for some $d'\in \Z$.
Substituting this into \eqref{eq:l1}, we see that $d'$ has to divide $b$. We write $b=d'b'$ with $b'\in \Z$.
We take $t=s/d'$ with some $s\in \Z$ which will be specified below.
We obtain
\begin{align*}
\beta_{s/d'}(P)
&=\frac{2c-as+b'd's x}{d'(bx-a)}\cdot \omega_W.
\end{align*}

Since $d'$ divides $b$, and since $(b,a)=1$, we have $(d',a)=1$. Therefore, we may chose $s\in \Z$ such that $2c\equiv as \pmod{d'}$. Then
$2c-as=ud'$ with some $u\in \Z$. We obtain that
\begin{align*}
\beta_{s/d'}(P)
&=\frac{u+b's x}{bx-a}\cdot \omega_W.
\end{align*}
Since $(a,b)=1$, this differential defines a rational section of $\Omega_{\calW/\Z}^1$ with no vertical poles.
This concludes the proof of the existence of $\underline{\alpha}(P)$.

If $\underline{\alpha}'(P)$ is another differential satisfying (i) and (ii), then $4(\underline{\alpha}(P)-\underline{\alpha}'(P))$
is a regular section of $\Omega_{\calW/\Z}^1$. Hence it is an integral multiple of $\omega_W$. This proves the uniqueness statement.
\end{proof}


We write $\underline\alpha_{\Delta,r}(f)$ for a differential on a minimal Weierstrass model of $E_{\Delta}$ corresponding to the point $P_{\Delta,r}(f)$ as in Lemma \ref{lem:wm}.

\begin{conjecture}
\label{con:ec}
Let the notation be as above.
Then
\[
\Delta c^{+}(\eps\Delta,r) - \frac{\eps c_E }{\Omega(E_\Delta)}\cdot \Re\int_{W_\Delta(\R)}\underline{\alpha}_{\Delta,r}(f)\in \frac{1}{4}\Z.
\]
\end{conjecture}

In Section \ref{subsect:ex} we present some numerical evidence for the conjecture, see Table \ref{table:3}.

\subsection{Examples}
\label{subsect:ex}

Here we illustrate  Corollary
\ref{cor:ec} in an example and present some numerical computations supporting Conjecture \ref{con:ec}.
We begin by briefly discussing how to compute the periods of differentials of the
third kind on elliptic curves as in Corollary \ref{cor:ec}.

Let $E:
y^2=4x^3-g_2 x -g_3$ be any elliptic curve over $\Q$.  If $Q=(x_0,y_0)\in E(\Q)$, we
consider the algebraic differential of the third kind $\beta(Q)$
with residue divisor  $\frac{1}{2}((Q)-(-Q))$ defined in \eqref{eq:beta}.
Let $L=\Z\mu+\Z\nu$ be the period lattice of $(E,\omega_E)$. We
assume that $E$ has real $2$-torsion (as in our example
\eqref{eq:defE} below).  Then we may chose $\mu$ and $\nu$ such that
$\mu\in \R$ and $\nu\in i\R$.  Let $\wp(z)$ be the Weierstrass
function corresponding to $L$.  The pullback of $\beta(Q)$ to the
complex uniformization $\C/L$ of $E$ is equal to
\[
 \frac{y_0}{\wp(z)-x_0} \cdot \frac{dz}{2}.
\]
Consequently, the  period $\Re\int_{E(\R)}\beta(Q)$ is equal to
\[
\Re \int_{u=0}^\mu  \frac{y_0}{\wp(u+iv)-x_0} \,\frac{du}{2}
\]
when $E(\R)$ is connected, and equal to twice this quantity when
$E(\R)$ has two components.  Here the integration is over the straight
line connecting $iv$ and $\mu+iv$ for a $v$ such that the line is
disjoint to $Q$ and $-Q$. Such integrals can be conveniently numerically computed using Maple.

We consider the elliptic curve
\begin{align}
\label{eq:defE}
E: y^2= 4x^3-4x+1
\end{align}
of conductor 37, which is the curve of smallest conductor with infinite Mordell-Weil group.
It has a modular parameterization by the modular curve $X_0(37)$, under which the N\'eron differential $\omega_E$ corresponds to the unique normalized Hecke eigenform $G\in S_2(\Gamma_0(37))$ which is invariant under the Fricke involution. So $\eps_G=+1$, and the Hecke $L$-function $L(G,s)$ has an odd functional equation.
The Manin constant of $E$ is equal to $1$ and the real period is given by
\[
\Omega(E)=5.98691729246391925966
\dots.
\]

Since the level $N$ is a prime, the space $S_{3/2,\bar\rho}$ is isomorphic to the space of scalar valued cusp forms of weight $3/2$ for $\Gamma_0(4\cdot 37)$ satisfying the Kohnen plus space condition. Analogously, $H_{1/2,\rho}$ is  isomorphic to the space of scalar valued harmonic Maass  forms of weight $1/2$ for $\Gamma_0(4\cdot 37)$ satisfying the Kohnen plus space condition.
In particular,
the coefficients of $f$,
the divisors $y_{\Delta,r}(f)$, and the differentials $\alpha_{\Delta,r}(f)$ only depend on $\Delta$ and not on the choice of $r$. Therefore we drop $r$ from the notation throughout this subsection.

We let $f$ be the unique harmonic Maass form of weight $1/2$ for $\Gamma_0(4\cdot 37)$ in the plus space with Fourier expansion
\[
f(\tau)=\sum_{n\gg -\infty} c^+(n) q^n + \sum_{n<0} c^-(n)\Gamma(1/2, 4\pi |n|
v) q^n,
\]
whose holomorphic part $f^+$ has the form
\[
f^+(\tau)=q^{-3}+O(q).
\]
Its image under $\xi$ is a non-zero multiple of the weight $3/2$ cusp form $g$ corresponding to $G$ under the Shimura correspondence.
For details on the numerical computation of the coefficients of $f$ we refer to \cite{BS}.

In Table \ref{table:2} we list for the first fundamental discriminants $\Delta$ (satisfying the Heegner hypothesis), the coefficient
$c^+(\Delta)$ of $f$, the Heegner point $P_\Delta(f)$ on the quadratic twist
\begin{align}
\label{eq:ed}
E_\Delta: y^2= 4x^3-4\Delta^2 x +\Delta^3,
\end{align}
and the difference
\[
\Delta c^+(\Delta) -\frac{1}{\Omega(E_\Delta)}\Re\left(\int_{E_\Delta(\R)} \alpha_\Delta(f)\right) .
\]
(For $\Delta= 77$ and $\Delta=101$ we omit some entries, since they would make the table too wide. These entries can easily be recovered from the corresponding entries of Table \ref{table:3}.)
The data was computed using Magma.
For instance, the first line of Table \ref{table:2} says that
\[
c^+(1)=\frac{1}{\Omega(E)}\Re\left(\int_{E(\R)} \frac{1}{2x}\, \frac{dx}{y}\right),
\]
and the numerical values confirm this up to 20 digits.

In Table \ref{table:3} we collect the data which is required for testing  Conjecture \ref{con:ec}.
A minimal Weierstrass model for the quadratic twist $E_\Delta$ is given by
\begin{align*}
W_\Delta: \begin{cases} y^2+y = x^3 -\Delta^2 x + (\Delta^3-1)/4, & \text{if $\Delta\equiv 1\pmod{4}$,}\\
y^2 = x^3 -\Delta^2 x + \Delta^3/4, & \text{if $\Delta\equiv 4\pmod{8}$,}\\
y^2 = x^3 -2^{-4}\Delta^2 x + 2^{-8}\Delta^3, & \text{if $\Delta\equiv 0\pmod{8}$.}
\end{cases}
\end{align*}
In the table we list for the first fundamental discriminants $\Delta$ (satisfying the Heegner hypothesis),
the Heegner point $P_\Delta(f)$ on $W_\Delta$, a $t\in \Q$ such that
the associated differential in \eqref{eq:l0} has no vertical poles,
and the difference
\[
\Delta c^+(\Delta) - \frac{1}{\Omega(E_\Delta)}\Re\left(\int_{W_\Delta(\R)} \underline\alpha_\Delta(f)\right).
\]
We find that the data is in accordance with the conjecture.
Note that the pullback of $\underline\alpha_\Delta(f)$ to $E_\Delta$ is equal to the differential
$\alpha_\Delta(f)+\frac{t}{2}\omega_{E_\Delta}$.


\begin{table}[h]
\caption{\label{table:2} Coefficients of $f$ and periods of $E_\Delta$}
\begin{tabular}{|r|c|c|c|c|c| }
\hline \rule[-3mm]{0mm}{8mm}
$\Delta$ &$c^+(\Delta)$& $P_\Delta(f)$ & difference\\
\hline 
\rule[-3mm]{0mm}{8mm}
$1$   &  $-0.28176 17849 89599 56879 \dots$
      &$(0,-1)$ & $0$\\
\rule[-3mm]{0mm}{8mm}
$12$  &  $-0.48852 72382 62012 25228 \dots$
      &$(1,-34)$& $-5$\\
\rule[-3mm]{0mm}{8mm}
$21$  &  $-0.17273 92572 32652 75652  \dots$
      &$(-\frac{335}{36},-\frac{16183}{108})$ & $-\frac{13}{6}$\\
\rule[-3mm]{0mm}{8mm}
$28$  &  $ \phantom{-}0.67819 39953 03947 79828 \dots$
      &$(-31,-2)$ & $19$\\
\rule[-3mm]{0mm}{8mm}
$33$  &  $ \phantom{-}0.56630 23201 59069 98168 \dots$
      &$(4, -137)$ & $20$\\
\rule[-3mm]{0mm}{8mm}
$37$  &  $ -0.91326 56137 46116 52958 \dots$
      &$(\frac{1009}{16},\frac{-26935}{32})$ &
      $-\frac{159}{4}$\\
\rule[-3mm]{0mm}{8mm}
$40$  &  $ \phantom{-}0.40098 50926 95436 37915 \dots$
      &$ (41,278)$ &
      $19$\\
\rule[-3mm]{0mm}{8mm}
$41$  &  $ \phantom{-}0.65637 49574 47572 31959 \dots$
      &$ (-\frac{344}{9} , \frac{8647}{27})$ &
      $\frac{77}{3}$\\
\rule[-3mm]{0mm}{8mm}
$44$  &  $ \phantom{-}0.96886 40443 45063 97321 \dots$
      &$(-\frac{46415}{1089}, -\frac{11674766}{35937}) $ &  $\frac{1445}{33}$\\
\rule[-3mm]{0mm}{8mm}
$53$  &  $  -0.56688 85256 82325 17859 \dots$
      &$ (-\frac{31839893674511}{880695910116}, -\frac{250032127988213200169}{413246299816000332})$ &
      $-\frac{26273369}{938454}$\\
\rule[-3mm]{0mm}{8mm}
$65$  &  $   -0.60328 07288 95214 77971 \dots$
      &$(\frac{27106}{225},-\frac{7720117}{3375}) $ &
      $-\frac{716}{15}$\\
\rule[-3mm]{0mm}{8mm}
$73$  &  $   \phantom{-}0.34874 71183 53624 08853 \dots$
      &$ (19,-107)$ &
      $26$\\
\rule[-3mm]{0mm}{8mm}
$77$ & $ \phantom{-}0.22699 13237 37052 54600 \dots $
&  & \\
\rule[-3mm]{0mm}{8mm}
$85$ & $-0.76894 61704 86762 72061 \dots$
& $( \frac{3649681}{24336},-\frac{5933836871}{1898208})$ & $-\frac{11651}{156}$\\
\rule[-3mm]{0mm}{8mm}
$101$ & $\phantom{-}0.24818 30869 47202 88319 \dots$
& &  \\
\hline
\end{tabular}

\end{table}

\begin{table}[h]
\caption{\label{table:3} The differentials $\underline\alpha_\Delta(f)$ and the associated periods}
\begin{tabular}{|r|c|c|c|c| }
\hline \rule[-3mm]{0mm}{8mm}
$\Delta$ & $P_\Delta(f)$ on $W_\Delta$ &  $t$
 & difference\\
\hline 
\rule[-3mm]{0mm}{8mm}
$1$ &   $(0,-1)$& $0$ & $0$\\
\rule[-3mm]{0mm}{8mm}
$12$&  $(1,-17)$& $0$  & $-5$\\
\rule[-3mm]{0mm}{8mm}
$21$ & $(-\frac{335}{36},-\frac{16291}{216})$ & $\frac{2}{3}$  & $-\frac{5}{2}$\\
\rule[-3mm]{0mm}{8mm}
$28$ &$(-31,-1)$ & $0$ &  $19$\\
\rule[-3mm]{0mm}{8mm}
$33$ &$(4, -69)$ & $0$ & $20$\\
\rule[-3mm]{0mm}{8mm}
$37$  & $(\frac{1009}{16},-\frac{26967}{64})$ & $\frac{1}{2}$ &
      $-40$\\
\rule[-3mm]{0mm}{8mm}
$40$ &
$ (\frac{41}{4},\frac{139}{8})$ & $0$ &
      $19$\\
\rule[-3mm]{0mm}{8mm}
$41$ &$ (-\frac{344}{9} , \frac{4310}{27})$ & $\frac{1}{3}$ &
      $\frac{51}{2}$\\
\rule[-3mm]{0mm}{8mm}
$44$ & $(-\frac{46415}{1089}, -\frac{5837383}{35937}) $ & $\frac{19}{33}$ &
      $\frac{87}{2}$\\
\rule[-3mm]{0mm}{8mm}
$53$ &$ (-\frac{31839893674511}{880695910116}, -\frac{250445374288029200501}{826492599632000664})$ &$\frac{3343}{469227}$&
      $-28$\\
\rule[-3mm]{0mm}{8mm}
$65$ & $(\frac{27106}{225},-\frac{3861746}{3375}) $ & $\frac{8}{15}$ &
      $-48$\\
\rule[-3mm]{0mm}{8mm}
$73$ & $ (19,-54)$ & $0$ &
      $26$\\
\rule[-3mm]{0mm}{8mm}
$77$ & $(\frac{1235881089099494174401}{82698967120806384516}, -\frac{128149977004435661308438102403131}{752055909346945729449269442936})$
& $\frac{4419608516}{4546948623}$ & $\frac{37}{2}$\\
\rule[-3mm]{0mm}{8mm}
$85$ & $(\frac{3649681}{24336}, -\frac{5933836871}{3796416})$  & $\frac{49}{78}$&  $-75$\\
\rule[-3mm]{0mm}{8mm}
$101$ & $(-\frac{4173521444186083063919}{290103403887032491044},
\frac{3127875148403162348360128767106003}
{4941163638674647451583112332072})$ & $\frac{7179562186}{8516211069}$ &  $\frac{43}{2}$\\
\hline
\end{tabular}

\end{table}


\begin{thebibliography}{BCDT}


\bibitem[Bo]{Bo1}
\emph{R. Borcherds}, Automorphic forms with singularities on
Grassmannians, Invent. Math. \textbf{132} (1998), 491--562.

\bibitem[BCDT]{BCDT} {\em C. Breuil, B. Conrad, F. Diamond, R. Taylor}, On the modularity of elliptic curves over $\Q$: wild 3-adic exercises,
J. Amer. Math. Soc. {\bf 14} (2001), 843--939.

\bibitem[BriO]{BriO} \emph{K. Bringmann and K. Ono},
Dyson's ranks and Maass forms, Ann. of Math. {\bf 171}  (2010),  419-–449.

\bibitem[Br]{Br} \emph{J. H. Bruinier}, Borcherds products on
  $\Orth(2,l)$ and Chern classes of Heegner divisors, Springer Lecture
  Notes in Mathematics {\bf 1780}, Springer-Verlag (2002).

\bibitem[BF]{BF} {\em J. H. Bruinier and J. Funke}, On two geometric theta lifts,
Duke Math. J. {\bf 125} (2004), 45--90.

\bibitem[BO]{BO} \emph{J. H. Bruinier and K. Ono}, Heegner divisors, L-functions and harmonic weak Maass forms, Annals of Math., accepted for publication.

\bibitem[BOR]{BOR} {\em J. H. Bruinier, K. Ono and R. Rhoades}, Differential operators for harmonic weak Maass forms and the vanishing of Hecke eigenvalues, Math. Ann. {\bf 342} (2008), 673--693.

\bibitem[BS]{BS} \emph{J. H. Bruinier and F. Str\"omberg},
Computation of harmonic Maass forms, preprint (2011).

\bibitem[Cr]{Cr} \emph{J. E. Cremona}, Algorithms for Modular Elliptic Curves, Cambridge University Press (1997).

\bibitem[Du]{Du} \emph{W. Duke},
Hyperbolic distribution problems and half-integral weight Maass forms,
Invent. Math. {\bf 92} (1988), 73–-90.

\bibitem[EZ]{EZ} \emph{M. Eichler and D. Zagier}, The Theory of Jacobi
  Forms, Progress in Math. {\bf 55}, Birkh\"auser (1985).

\bibitem[Ed]{Ed} \emph{B. Edixhoven}, On the Manin constants of modular elliptic curves. In: Arithmetic Algebraic Geometry, Texel 1989 (eds.: G. van der Geer, F Oort, and J. Steenbrink), 25--39, Progr. Math. {\bf 89}, Boston, Birh\"auser (1991).

\bibitem[GKZ]{GKZ} \emph{B. Gross, W. Kohnen, and D. Zagier}, Heegner
  points and derivatives of $L$-series. II.  Math. Ann.  {\bf 278}
  (1987), 497--562.

\bibitem[GH]{Gr} \emph{P. A. Griffiths and J. Harris}, Principles of Algebraic Geometry,
John Wiley  (1978).

\bibitem[GZ]{GZ} {\em B. Gross and D. Zagier}, Heegner points and derivatives of
L-series, Invent. Math. {\bf 84} (1986), 225--320.

\bibitem[KS]{KS} \emph{S. Katok and P. Sarnak},
Heegner points, cycles and Maass forms, Israel J. Math. {\bf 84}
(1993), 193--227.


\bibitem[KZ]{KZ} \emph{W. Kohnen and D. Zagier},
Values of $L$-series of modular forms at the center of the critical
strip.  Invent. Math.  {\bf 64}  (1981), no. 2, 175--198.

\bibitem[KonZ]{KonZ} \emph{M. Kontsevich and D. Zagier}, Periods.
Mathematics unlimited---2001 and beyond,  771--808, Springer, Berlin (2001).

\bibitem[On]{On} {\emph K. Ono}, Unearthing the visions of a master:
harmonic Maass forms and number theory,  Current developments in mathematics, 2008,  347-–454,
Int. Press, Somerville (2009).

\bibitem[Sch]{Sch} {\em A. J. Scholl}, Fourier coefficients of Eisenstein series
 on non-congruence subgroups, Math. Proc. Camb. Phil. Soc.
  {\bf 99} (1986), 11--17.

\bibitem[Sh]{Sh}
\emph{G. Shimura}, On modular forms of half integral weight. Ann. of
Math. (2) {\bf 97} (1973), 440--481.

\bibitem[Si]{Si} {\em J. H. Silverman}, The arithmetic of elliptic curves.
Graduate Texts in Mathematics {\bf 106}, Springer, Dordrecht (2009).
\bibitem[Sk1]{Sk1} \emph{N.-P. Skoruppa}, Developments in the theory
of Jacobi forms. In: Proceedings of the conference on automorphic
funtions and their applications, Chabarovsk (eds.: N. Kuznetsov and
V. Bykovsky), The USSR Academy of Science (1990), 167--185. (see
also MPI-preprint 89-40, Bonn (1989).)


\bibitem[Sk2]{Sk2} \emph{N.-P. Skoruppa}, Explicit formulas for the
  Fourier coefficients of Jacobi and elliptic modular forms.  Invent.
  Math.  {\bf 102} (1990), 501--520.

\bibitem[SZ]{SZ} \emph{N.-P. Skoruppa and D. Zagier}, Jacobi forms
and a certain space of modular forms,  Invent. Math.  {\bf 94}
(1988), 113--146.

\bibitem[W]{W} \emph{M. Waldschmidt}, Nombers transcendents et
groupes alg\'ebraiques, Ast\'erisque {\bf 69--70} (1979).

\bibitem[Wa]{Wa}
\emph{J.-L. Waldspurger}, Sur les coefficients de Fourier des formes
modulaires de poids demi-entier,
J. Math. Pures Appl. (9) {\bf 60}  (1981), no. 4, 375--484.

\bibitem[Wi]{Wi} {\em A. Wiles},
Modular elliptic curves and Fermat's last theorem,
Ann. of Math. {\bf 141} (1995), 443--551.

\bibitem[Za]{Za} \emph{D. Zagier}, Ramanujan's mock theta functions
and their applications [d'apr\`es Zwegers and
Bringmann-Ono], S\'eminaire Bourbaki 60\'eme ann\'ee, 2006-2007,
no. 986.

\bibitem[Zw1]{Z1} \emph{S. P. Zwegers}, {Mock $\vartheta$-functions
 and real analytic modular forms}, $q$-series with applications to
 combinatorics, number theory, and physics (Ed. B. C. Berndt and K.
 Ono), Contemp. Math. \textbf{291}, Amer. Math. Soc., (2001),
 269--277.

 \bibitem[Zw2]{Z2} \emph{S. P. Zwegers}, {Mock theta functions},
 Ph.D. Thesis, Universiteit Utrecht, 2002.


\end{thebibliography}
\end{document}